\documentclass[11pt]{article}

\usepackage[a4paper,margin=1in]{geometry}
\usepackage{amsmath,amssymb}
\usepackage{url}

\title{A Regularised Wallis Hierarchy}
\author{
  S. R. Holcombe\\
  MBS, University of Melbourne\\
}
\date{}

\begin{document}

\maketitle

\begin{abstract}
A hierarchy of regularised Wallis products is introduced by raising the reciprocal Wallis
factor
\[
  1-\frac1{n^2}
\]
to the polynomial weight \(n^m\), \(m=0,1,2,\ldots\).  For each \(m\), a minimal exponential
counterterm is chosen by cancelling precisely the non-summable terms in the logarithmic
expansion.  This gives a convergent product \(P_m\) the logarithm of which is an explicit
zeta-function tail.  The first non-trivial examples are
  \[
  \prod_{n=2}^{\infty}
  e^{1/n}
  \left(1-\frac1{n^2}\right)^n
  =
  \frac{e^\gamma}{2},
  \qquad
  \prod_{n=2}^{\infty}
  e\left(1-\frac1{n^2}\right)^{n^2}
  =
  \frac{\pi}{e^{3/2}}.
  \]
The even branch has a finite closed form involving \(\pi\), harmonic numbers, and odd zeta
values.  The odd branch reduces to finite logarithmic gamma moments, and hence to constants
involving \(\gamma\), logarithms, odd zeta values, and derivatives of the zeta function at
positive even integers.  The same subtraction rule also gives a two-factor extension
involving the companion factor \(1+1/n^2\).  Finally, the associated \(x\)-dependent
products factor into one-sided canonical products, giving a direct connection with
Kurokawa's multiple sine functions: the even Wallis branch is obtained from odd multiple
sine functions, while the odd branch appears as a symmetric companion to the even multiple
sine case.
\end{abstract}

\section{Introduction}

Wallis' product is the classical identity \cite{wallis}
\[
  \frac{\pi}{2} = \prod_{n=1}^{\infty} \frac{4n^2}{4n^2-1}.
\]
The reciprocal factor gives the elementary telescoping product
\[
  \prod_{n=2}^{\infty} \left(1-\frac1{n^2}\right) = \frac12.
\]
The present note starts from this reciprocal factor and asks what happens when its exponent
is no longer constant.  The replacement
\[
  \left(1-\frac1{n^2}\right) \quad\longmapsto\quad \left(1-\frac1{n^2}\right)^{n^m}
\]
is simple, but it immediately destroys ordinary convergence for most \(m\).  The logarithm
of each factor has an explicit asymptotic expansion, so the divergent part can be identified
and removed term by term.

A previous short note gave the product representation \cite{holcombe}
\[
  \prod_{n=2}^{\infty} e\left(1-\frac1{n^2}\right)^{n^2} = \frac{\pi}{e^{3/2}},
\]
which was subsequently recovered by Allouche \cite{allouche} from Kurokawa's triple sine
function.  In the notation used below it is the \(m=2\) member of the hierarchy.

The purpose of the present note is to place this product in a simple indexed family.  The
construction extends the reciprocal Wallis product by replacing the exponent \(1\) with
\(n^m\).  This introduces divergent logarithmic terms, so each product is regularised by
subtracting exactly those terms and no summable terms.  Equivalently, the logarithmic
divergence is removed from
\[
  n^m\log\left(1-\frac1{n^2}\right)
\]
by a minimal exponential counterterm.

Wallis-type products and zeta-regularised products are classical objects, and there are many
generalisations of Wallis' formula.  Nearby Wallis-type generalisations include the work of
Cai, Hu and Kim \cite{caihukim}, and Farrell \cite{farrell}.  The connection with Kurokawa's
multiple trigonometric functions is also relevant here: the first non-trivial even member is
the triple-sine case identified by Allouche, and the even branch below is naturally related
to Kurokawa's odd multiple sine functions
\cite{allouche,kurokawawakayama1,kurokawawakayama2}.  The present formulation is different
in a narrower sense.  It applies a fixed minimal-subtraction rule to the same Wallis factor
for every \(m\), and this produces a regularised Wallis hierarchy whose logarithms are
explicit zeta tails.

\section{The hierarchy and its zeta tails}

For \(m=0,1,2,\ldots\), put
\[
  J_m = \left\lfloor \frac{m+1}{2}\right\rfloor .
\]
Define
\[
  P_m = \prod_{n=2}^{\infty} \exp\left( \sum_{j=1}^{J_m} \frac{n^{m-2j}}{j} \right) \left(1-\frac1{n^2}\right)^{n^m}.
\]
This is the regularised Wallis hierarchy.  Each product is built from the same Wallis
factor, but with increasing polynomial weights \(n^m\).  The exponential factor is the
minimal counterterm: it cancels exactly the terms in the logarithmic expansion whose sums
over \(n\) are divergent.

Indeed,
\[
  \log\left(1-\frac1{n^2}\right) = -\sum_{j=1}^{\infty} \frac{1}{j n^{2j}},
\]
and hence
\[
  n^m\log\left(1-\frac1{n^2}\right) = -\sum_{j=1}^{\infty} \frac{n^{m-2j}}{j}.
\]
The term \(n^{m-2j}\) is summable over \(n\) when
\[
  m-2j<-1.
\]
Consequently the terms requiring subtraction are exactly those with
\[
  j\leq \frac{m+1}{2}.
\]
This gives
\[
  J_m = \left\lfloor \frac{m+1}{2}\right\rfloor .
\]

For the finite product
\[
  P_{m,N} = \prod_{n=2}^{N} \exp\left( \sum_{j=1}^{J_m} \frac{n^{m-2j}}{j} \right) \left(1-\frac1{n^2}\right)^{n^m},
\]
the cancellation gives
\[
  \log P_{m,N} = -\sum_{n=2}^{N} \sum_{j=J_m+1}^{\infty} \frac{n^{m-2j}}{j}.
\]
After taking \(N\to\infty\), the remaining double series is absolutely convergent.  The
order of summation may therefore be interchanged:
\[
  \log P_m = -\sum_{j=J_m+1}^{\infty} \frac1j \sum_{n=2}^{\infty} \frac1{n^{2j-m}}.
\]
Hence
\[
  \log P_m = -\sum_{j=J_m+1}^{\infty} \frac{\zeta(2j-m)-1}{j}.
\]

The parity of \(m\) separates the hierarchy into two natural branches.  If \(m=2a\),
\(a\geq1\), then
\[
  \log P_{2a} = -\sum_{r=1}^{\infty} \frac{\zeta(2r)-1}{r+a}.
\]
If \(m=2a+1\), \(a\geq0\), then
\[
  \log P_{2a+1} = -\sum_{r=1}^{\infty} \frac{\zeta(2r+1)-1}{r+a+1}.
\]

\section{Evaluation of the two branches}

For the even branch define
\[
  T_{2a} = \sum_{r=1}^{\infty} \frac{\zeta(2r)-1}{r+a}, \qquad \log P_{2a}=-T_{2a}.
\]
Euler's product for the sine function gives
\[
  \sum_{r=1}^{\infty} \zeta(2r)x^{2r} = \frac{1-\pi x\cot(\pi x)}{2},
\]
and hence
\[
  \sum_{r=1}^{\infty} \bigl(\zeta(2r)-1\bigr)x^{2r} = \frac{1-\pi x\cot(\pi x)}{2} - \frac{x^2}{1-x^2}.
\]
Using
\[
  \frac1{r+a} = 2\int_0^1 x^{2r+2a-1}\,dx,
\]
one obtains
\[
  T_{2a}
  =
  2\int_0^1
  x^{2a-1}
  \left[
    \frac{1-\pi x\cot(\pi x)}{2}
    -
    \frac{x^2}{1-x^2}
    \right]dx.
\]
The integrand may be written as a logarithmic derivative, since
\[
  \pi\cot(\pi x)+\frac{2x}{1-x^2} = \frac{d}{dx} \log\left(\frac{\sin(\pi x)}{1-x^2}\right).
\]
An integration by parts then gives
\[
  T_{2a}
  =
  \frac1{2a}
  -
  \log\left(\frac{\pi}{2}\right)
  +
  2a
  \int_0^1
  x^{2a-1}
  \log\left(\frac{\sin(\pi x)}{1-x^2}\right)dx.
\]
Writing
\[
  I_a= \int_0^1 x^{2a-1}\log(\sin(\pi x))\,dx, \qquad K_a= \int_0^1 x^{2a-1}\log(1-x^2)\,dx,
\]
this becomes
\[
  T_{2a}
  =
  \frac1{2a}
  -
  \log\left(\frac{\pi}{2}\right)
  +
  2aI_a
  -
  2aK_a.
\]
The elementary moment is
\[
  K_a=-\frac{H_a}{2a}.
\]
The log-sine moment is
\[
  I_a
  =
  -\frac{\log2}{2a}
  -
  \sum_{\ell=1}^{a-1}
  (-1)^{\ell+1}
  \frac{(2a-1)!}{(2a-2\ell)!}
  \frac{\zeta(2\ell+1)}{(2\pi)^{2\ell}}.
\]
Substitution gives the even branch in closed form:
\[
  \log P_{2a}
  =
  \log\pi
  -
  H_a
  -
  \frac1{2a}
  +
  2a
  \sum_{\ell=1}^{a-1}
  (-1)^{\ell+1}
  \frac{(2a-1)!}{(2a-2\ell)!}
  \frac{\zeta(2\ell+1)}{(2\pi)^{2\ell}}.
\]
Equivalently,
\[
  P_{2a}
  =
  \pi
  \exp\left[
    -
    H_a
    -
    \frac1{2a}
    +
    2a
    \sum_{\ell=1}^{a-1}
    (-1)^{\ell+1}
    \frac{(2a-1)!}{(2a-2\ell)!}
    \frac{\zeta(2\ell+1)}{(2\pi)^{2\ell}}
    \right].
\]

For the odd branch define
\[
  T_{2a+1} = \sum_{r=1}^{\infty} \frac{\zeta(2r+1)-1}{r+a+1}, \qquad \log P_{2a+1}=-T_{2a+1}.
\]
The standard expansion
\[
  \log\Gamma(1+x) = -\gamma x + \sum_{k=2}^{\infty} \frac{(-1)^k\zeta(k)}{k}x^k
\]
gives
\[
  \sum_{r=1}^{\infty} \zeta(2r+1)x^{2r} = -\frac12 \left[ \psi(1+x)+\psi(1-x)+2\gamma \right].
\]
Therefore
\[
  \sum_{r=1}^{\infty}
  \bigl(\zeta(2r+1)-1\bigr)x^{2r}
  =
  -\frac12
  \left[
    \psi(1+x)+\psi(1-x)+2\gamma
    \right]
  -
  \frac{x^2}{1-x^2}.
\]
Using
\[
  \frac1{r+a+1} = 2\int_0^1 x^{2r+2a+1}\,dx,
\]
the odd tail becomes
\[
  \log P_{2a+1}
  =
  \int_0^1
  x^{2a+1}
  \left[
    \psi(1+x)+\psi(1-x)+2\gamma
    +
    \frac{2x^2}{1-x^2}
    \right]dx.
\]
Although the last two terms in the bracket are separately singular at \(x=1\), their
singular parts cancel.  The pole of \(\psi(1-x)\) is cancelled by the pole of
\(2x^2/(1-x^2)\), leaving a finite integrand at the endpoint.  Thus the displayed integral
is an ordinary convergent integral, with the cancellation kept inside the bracket.

For each fixed \(a\), the odd-branch integral can be evaluated in closed finite form.
Define
\[
  M_k = \int_0^1 x^k\left[\psi(1+x)+\gamma\right]dx, \qquad k=0,1,2,\ldots .
\]
Then
\[
  \log P_{2a+1}
  =
  \sum_{k=0}^{2a}
  (-1)^k
  \binom{2a+1}{k}
  M_k
  +
  H_{2a+1}
  -
  H_a
  -
  \log2
  -
  \frac1{a+1}.
\]
Here
\[
  M_0=\gamma,
\]
and, for \(k\geq1\),
\[
  M_k = \frac{\gamma}{k+1} + \frac1k - kL_{k-1},
\]
where
\[
  L_p = \int_0^1x^p\log\Gamma(x)\,dx.
\]
Kummer's Fourier expansion for \(\log\Gamma(x)\) \cite{whittakerwatson}, gives
\[
  \begin{aligned}
    L_p
     & =
    (\gamma+\log2)
    \left[
      \frac1{2(p+1)}
      -
      \frac1{p+2}
      \right]
    +
    \log\pi
    \left[
      \frac1{p+1}
      -
      \frac1{p+2}
      \right]
    \\
     & \quad
    +
    \frac{\log2}{2(p+1)}
    +
    \frac12
    \sum_{\ell=1}^{\lfloor p/2\rfloor}
    (-1)^{\ell+1}
    \frac{p!}{(p-2\ell+1)!}
    \frac{\zeta(2\ell+1)}{(2\pi)^{2\ell}}
    \\
     & \quad
    +
    \frac1\pi
    \sum_{\ell=0}^{\lfloor (p-1)/2\rfloor}
    (-1)^\ell
    \frac{p!}{(p-2\ell)!}
    \frac{\zeta'(2\ell+2)}{(2\pi)^{2\ell+1}}.
  \end{aligned}
\]
Thus the odd branch also reduces to finite expressions.  For each fixed \(a\), the product
\(P_{2a+1}\) is expressed in terms of \(\gamma\), \(\log2\), \(\log\pi\), odd zeta values,
and derivatives of \(\zeta(s)\) at positive even integers.

Combining the two branches, the regularised Wallis hierarchy for even indices \(m=2a\)
\(a\geq1\), is
\[
  \prod_{n=2}^{\infty}
  \exp\left(
  \sum_{j=1}^{a}
  \frac{n^{2a-2j}}{j}
  \right)
  \left(1-\frac1{n^2}\right)^{n^{2a}}
  =
  \pi
  \exp\left[
    -
    H_a
    -
    \frac1{2a}
    +
    2a
    \sum_{\ell=1}^{a-1}
    (-1)^{\ell+1}
    \frac{(2a-1)!}{(2a-2\ell)!}
    \frac{\zeta(2\ell+1)}{(2\pi)^{2\ell}}
    \right].
\]
For odd indices, \(m=2a+1\), \(a\geq0\),
\[
  \prod_{n=2}^{\infty}
  \exp\left(
  \sum_{j=1}^{a+1}
  \frac{n^{2a+1-2j}}{j}
  \right)
  \left(1-\frac1{n^2}\right)^{n^{2a+1}}
  =
  \exp\left[
    \sum_{k=0}^{2a}
    (-1)^k
    \binom{2a+1}{k}
    M_k
    +
    H_{2a+1}
    -
    H_a
    -
    \log2
    -
    \frac1{a+1}
    \right],
\]
where \(M_0=\gamma\), and, for \(k\geq1\),
\[
  M_k
  =
  \frac{\gamma}{k+1}
  +
  \frac1k
  -
  kL_{k-1},
  \qquad
  L_p
  =
  \int_0^1x^p\log\Gamma(x)\,dx,
\]
with \(L_p\) given explicitly above.

\section{Hyperbolic companion}

The construction can be extended in several directions.  Here however the focus is on
keeping the same powers \(n^m\) and the same minimal-subtraction rule, but to allow two
Wallis-like factors.  The reciprocal Wallis factor
\[
  1-\frac1{n^2}
\]
is now accompanied by its hyperbolic companion
\[
  1+\frac1{n^2}.
\]
Thus the two-factor product is built from
\[
  1-\frac1{n^2}, \qquad 1+\frac1{n^2}.
\]

Let \(A\) and \(B\) be fixed parameters.  For \(m=0,1,2,\ldots\), put
\[
  J_m = \left\lfloor\frac{m+1}{2}\right\rfloor .
\]
Define
\[
\begin{aligned}
  P_m(A,B)
  &=
  \prod_{n=2}^{\infty}
  \exp\left(
  \sum_{j=1}^{J_m}
  \frac{A+(-1)^jB}{j}
  n^{m-2j}
  \right)
  \\
  &\quad\times
  \left[
  \left(1-\frac1{n^2}\right)^A
  \left(1+\frac1{n^2}\right)^B
  \right]^{n^m}.
\end{aligned}
\]
The original hierarchy is the special case
\[
  P_m=P_m(1,0).
\]

The coefficient \(A+(-1)^jB\) comes from
\[
  A\log(1-u)+B\log(1+u) = - \sum_{j=1}^{\infty} \frac{A+(-1)^jB}{j}u^j, \qquad |u|<1.
\]
With \(u=1/n^2\), multiplication by \(n^m\) gives
\[
n^m
\left[
A\log\left(1-\frac1{n^2}\right)
+
B\log\left(1+\frac1{n^2}\right)
\right]
=
-
\sum_{j=1}^{\infty}
\frac{A+(-1)^jB}{j}
n^{m-2j}.
\]
The non-summable powers are again those with
\[
  m-2j\geq -1.
\]
Thus the same cutoff \(J_m\) removes the divergent part of the logarithm.  For special
choices of \(A\) and \(B\), some of the coefficients may vanish, but the same expression
gives a uniform two-parameter form.

After cancellation, the remaining logarithm is absolutely summable.  The zeta-tail formula
becomes
\[
  \log P_m(A,B) = - \sum_{j=J_m+1}^{\infty} \frac{A+(-1)^jB}{j} \left(\zeta(2j-m)-1\right).
\]
The dependence on \(A\) and \(B\) is therefore linear at the logarithmic level.

The hyperbolic companion is
\[
  P_m^+ = P_m(0,1).
\]
Explicitly,
\[
P_m^+
=
\prod_{n=2}^{\infty}
\exp\left(
\sum_{j=1}^{J_m}
\frac{(-1)^j}{j}
n^{m-2j}
\right)
\left(1+\frac1{n^2}\right)^{n^m},
\]
and
\[
  \log P_m^+ = - \sum_{j=J_m+1}^{\infty} \frac{(-1)^j}{j} \left(\zeta(2j-m)-1\right).
\]
Consequently
\[
  \log P_m(A,B) = A\log P_m+B\log P_m^+,
\]
or
\[
  P_m(A,B)=P_m^A(P_m^+)^B.
\]

For even indices, put \(m=2a\), with \(a\geq1\).  Then
\[
  \log P_{2a}(A,B) = A\log P_{2a} + B\log P_{2a}^+.
\]
The first term is the even branch obtained previously:
\[
\log P_{2a}
=
\log\pi
-
H_a
-
\frac1{2a}
+
2a
\sum_{\ell=1}^{a-1}
(-1)^{\ell+1}
\frac{(2a-1)!}{(2a-2\ell)!}
\frac{\zeta(2\ell+1)}{(2\pi)^{2\ell}}.
\]

Using similar methods as before, the corresponding hyperbolic factor can be shown to be
\[
\begin{aligned}
  \log P_{2a}^+
  &=
  -\log2
  -
  \sum_{j=1}^{a}\frac{(-1)^j}{j}
  +
  (-1)^a
  \left[
  \frac{\pi}{2a+1}
  -
  \frac1{2a}
  \right.
  \\
  &\qquad\left.
  +
  \frac{(2a)!}{(2\pi)^{2a}}
  \left\{
  \zeta(2a+1)
  -
  \sum_{s=0}^{2a}
  \frac{(2\pi)^s}{s!}
  \operatorname{Li}_{2a+1-s}(e^{-2\pi})
  \right\}
  \right].
\end{aligned}
\]
Thus the even part of the two-factor hierarchy is

\[
\begin{aligned}
  P_{2a}(A,B)
  =
  \exp\Bigg\{&
  A
  \left[
  \log\pi
  -
  H_a
  -
  \frac1{2a}
  +
  2a
  \sum_{\ell=1}^{a-1}
  (-1)^{\ell+1}
  \frac{(2a-1)!}{(2a-2\ell)!}
  \frac{\zeta(2\ell+1)}{(2\pi)^{2\ell}}
  \right]
  \\
  &+
  B
  \left[
  -\log2
  -
  \sum_{j=1}^{a}\frac{(-1)^j}{j}
  +
  (-1)^a
  \left(
  \frac{\pi}{2a+1}
  -
  \frac1{2a}
  \right.
  \right.
  \\
  &\qquad\left.\left.
  +
  \frac{(2a)!}{(2\pi)^{2a}}
  \left\{
  \zeta(2a+1)
  -
  \sum_{s=0}^{2a}
  \frac{(2\pi)^s}{s!}
  \operatorname{Li}_{2a+1-s}(e^{-2\pi})
  \right\}
  \right)
  \right]
  \Bigg\}.
\end{aligned}
\]

\section{One-sided factors and Kurokawa multiple sine functions}

The products above are symmetric in the two factors \(1-x/n\) and \(1+x/n\).  It is useful
to separate this symmetric product into one-sided canonical factors.  This factorisation
gives a convenient framework for comparing nearby infinite products built from one-sided
factors, and also explains how the ordinary Wallis products and their hyperbolic companions
sit inside the same Kurokawa multiple-sine framework.

Define the \(x\)-dependent product
\[
F_m(x)
=
\prod_{n=1}^{\infty}
\exp\left(
\sum_{j=1}^{J_m}
\frac{x^{2j}n^{m-2j}}{j}
\right)
\left(1-\frac{x^2}{n^2}\right)^{n^m}.
\]
This is the regularised Wallis product with \(1\) replaced by \(x^2\).

The factor with \(n=1\) is
\[
  \exp\left( \sum_{j=1}^{J_m} \frac{x^{2j}}{j} \right) (1-x^2).
\]
Thus \(F_m(x)\) has a simple zero at \(x=1\).  Writing
\[
  F_m(x) = \exp\left( \sum_{j=1}^{J_m} \frac{x^{2j}}{j} \right) (1-x^2)G_m(x),
\]
where
\[
G_m(x)
=
\prod_{n=2}^{\infty}
\exp\left(
\sum_{j=1}^{J_m}
\frac{x^{2j}n^{m-2j}}{j}
\right)
\left(1-\frac{x^2}{n^2}\right)^{n^m},
\]
then
\[
  G_m(1)=P_m.
\]
Differentiating at the simple zero gives
\[
  F_m'(1) = \exp\left( \sum_{j=1}^{J_m}\frac1j \right) \left.\frac{d}{dx}(1-x^2)\right|_{x=1} G_m(1).
\]
Since
\[
  \sum_{j=1}^{J_m}\frac1j=H_{J_m}, \qquad \left.\frac{d}{dx}(1-x^2)\right|_{x=1}=-2,
\]
it follows that
\[
  P_m = -\frac{F_m'(1)}{2e^{H_{J_m}}}.
\]

The hyperbolic companion is obtained from the same product by putting \(x=i\).  Indeed,
\[
F_m(i)
=
\prod_{n=1}^{\infty}
\exp\left(
\sum_{j=1}^{J_m}
\frac{(-1)^j n^{m-2j}}{j}
\right)
\left(1+\frac1{n^2}\right)^{n^m}.
\]
The \(n=1\) factor is
\[
  2\exp\left( \sum_{j=1}^{J_m}\frac{(-1)^j}{j} \right).
\]
Therefore, if
\[
  E_m = \sum_{j=1}^{J_m}\frac{(-1)^j}{j},
\]
then
\[
  P_m^+ = \frac{F_m(i)}{2e^{E_m}}.
\]
Thus the ordinary Wallis product \(P_m\) is obtained from the behaviour of \(F_m(x)\) at the
zero \(x=1\), while the hyperbolic companion \(P_m^+\) is obtained from the value of the
same product at the imaginary point \(x=i\).

Next define the one-sided canonical product
\[
\Phi_m(z)
=
\prod_{n=1}^{\infty}
\exp\left(
\sum_{j=1}^{m+1}
\frac{(-1)^jz^j n^{m-j}}{j}
\right)
\left(1+\frac zn\right)^{n^m}.
\]
This is a one-sided canonical product associated with the same subtraction principle.
Indeed,
\[
  n^m\log\left(1+\frac zn\right) = \sum_{j=1}^{\infty} (-1)^{j+1} \frac{z^j n^{m-j}}{j},
\]
so the exponential factor in \(\Phi_m(z)\) cancels the first \(m+1\) terms of this
logarithmic expansion.

Now form the symmetric product
\[
  \Phi_m(x)\Phi_m(-x).
\]
The algebraic factors give
\[
  \left(1+\frac xn\right)^{n^m} \left(1-\frac xn\right)^{n^m} = \left(1-\frac{x^2}{n^2}\right)^{n^m}.
\]
The exponential counterterms add to
\[
  \sum_{j=1}^{m+1} \frac{\left[(-1)^jx^j+(-1)^j(-x)^j\right]n^{m-j}}{j}.
\]
The odd powers cancel.  If \(j=2q\), the surviving term is
\[
  \frac{x^{2q}n^{m-2q}}{q}.
\]
The allowed values of \(q\) are
\[
  1\leq q\leq \left\lfloor\frac{m+1}{2}\right\rfloor = J_m.
\]
Therefore
\[
  \Phi_m(x)\Phi_m(-x)=F_m(x).
\]

Kurokawa's multiple sine functions are defined from the canonical factors
\[
  \mathcal P_r(u) = (1-u) \exp\left( u+\frac{u^2}{2}+\cdots+\frac{u^r}{r} \right).
\]
With the normalisation used by Kurokawa and Wakayama, one has for \(r\geq2\),
\[
S_r(z)
=
\exp\left(\frac{z^{r-1}}{r-1}\right)
\prod_{n=1}^{\infty}
\left[
\mathcal P_r\!\left(\frac zn\right)
\mathcal P_r\!\left(-\frac zn\right)^{(-1)^{r-1}}
\right]^{n^{r-1}} .
\]
See \cite{allouche,kurokawawakayama2}.  Taking \(r=m+1\), this gives
\[
S_{m+1}(z)
=
\exp\left(\frac{z^m}{m}\right)
\prod_{n=1}^{\infty}
\left[
\mathcal P_{m+1}\!\left(\frac zn\right)
\mathcal P_{m+1}\!\left(-\frac zn\right)^{(-1)^m}
\right]^{n^m}.
\]
Since
\[
  \Phi_m(z) = \prod_{n=1}^{\infty} \left[ \mathcal P_{m+1}\!\left(-\frac zn\right) \right]^{n^m},
\]
Kurokawa's product may be written as
\[
  S_{m+1}(z) = e^{z^m/m} \Phi_m(-z)\Phi_m(z)^{(-1)^m}, \qquad m\geq1.
\]

For even \(m=2a\), with \(a\geq1\), this becomes
\[
  S_{2a+1}(z) = e^{z^{2a}/(2a)} \Phi_{2a}(-z)\Phi_{2a}(z).
\]
Using
\[
  \Phi_{2a}(z)\Phi_{2a}(-z)=F_{2a}(z),
\]
yields
\[
  F_{2a}(z) = e^{-z^{2a}/(2a)}S_{2a+1}(z).
\]
This identity gives both the ordinary even branch and its hyperbolic companion.

At \(z=1\), \(F_{2a}(1)=0\), so
\[
  S_{2a+1}(1)=0.
\]
Differentiating
\[
  F_{2a}(z) = e^{-z^{2a}/(2a)}S_{2a+1}(z)
\]
and evaluating at \(z=1\) gives
\[
  F_{2a}'(1) = e^{-1/(2a)}S_{2a+1}'(1).
\]
Since
\[
  P_{2a} = -\frac{F_{2a}'(1)}{2e^{H_a}},
\]
then
\[
  S_{2a+1}'(1) = -2\exp\left(H_a+\frac1{2a}\right)P_{2a}.
\]
Substituting the closed form for \(P_{2a}\) gives
\[
S_{2a+1}'(1)
=
-2\pi
\exp\left[
2a
\sum_{\ell=1}^{a-1}
(-1)^{\ell+1}
\frac{(2a-1)!}{(2a-2\ell)!}
\frac{\zeta(2\ell+1)}{(2\pi)^{2\ell}}
\right].
\]
Note that for \(a=1\)
\[
  S_3'(1)=-2\pi,
\]
and hence
\[
  P_2 = -\frac{S_3'(1)} {2\exp\left(H_1+\frac12\right)} = \frac{\pi}{e^{3/2}}.
\]
This is the case identified by Allouche through Kurokawa's triple sine function
\cite{allouche}.

At the imaginary point \(z=i\), the same identity gives
\[
  F_{2a}(i) = e^{-i^{2a}/(2a)}S_{2a+1}(i) = e^{-(-1)^a/(2a)}S_{2a+1}(i).
\]
Since
\[
  P_{2a}^+ = \frac{F_{2a}(i)} {2\exp\left(\sum_{j=1}^{a}(-1)^j/j\right)},
\]
then
\[
S_{2a+1}(i)
=
2
\exp\left(
\sum_{j=1}^{a}\frac{(-1)^j}{j}
+
\frac{(-1)^a}{2a}
\right)
P_{2a}^+.
\]
Thus the even hyperbolic companion is the value of the corresponding Kurokawa odd multiple
sine at the imaginary argument \(i\), after removing the elementary \(n=1\) factor.

For odd \(m=2a+1\), Kurokawa's product gives instead
\[
  S_{2a+2}(z) = e^{z^{2a+1}/(2a+1)} \frac{\Phi_{2a+1}(-z)}{\Phi_{2a+1}(z)}.
\]
The Kurokawa multiple sine of even order therefore involves the quotient of the two
one-sided canonical factors.  The Wallis hierarchy and its hyperbolic companion instead
involve the symmetric products
\[
  F_{2a+1}(1) = \Phi_{2a+1}(1)\Phi_{2a+1}(-1),
\]
and
\[
  F_{2a+1}(i) = \Phi_{2a+1}(i)\Phi_{2a+1}(-i).
\]
Thus the odd branch is naturally viewed as the symmetric companion to Kurokawa's even
multiple sine, rather than as a Kurokawa multiple sine directly.

\section{Examples}

The first few members of the hierarchy are as follows.  The zeroth member is the telescoping
Wallis factor,
\[
  P_0 = \prod_{n=2}^{\infty} \left(1-\frac1{n^2}\right) = \frac12.
\]
The first odd member is
\[
  P_1 = \prod_{n=2}^{\infty} e^{\frac1n} \left(1-\frac1{n^2}\right)^n = \frac{e^\gamma}{2}.
\]
The first non-trivial even member is
\[
  P_2 = \prod_{n=2}^{\infty} e \left(1-\frac1{n^2}\right)^{n^2} = \frac{\pi}{e^{3/2}}.
\]
The next odd member is
\[
  P_3
  =
  \prod_{n=2}^{\infty}
  e^{n+\frac{1}{2n}}
  \left(1-\frac1{n^2}\right)^{n^3}
  =
  \sqrt{\frac{\pi}{2}}\,
  \exp\left(
  \gamma-\frac76-\frac{3\zeta'(2)}{\pi^2}
  \right).
\]
The next even member is
\[
  P_4
  =
  \prod_{n=2}^{\infty}
  e^{n^2+\frac12}
  \left(1-\frac1{n^2}\right)^{n^4}
  =
  \pi
  \exp\left(
  -\frac74
  +
  \frac{3\zeta(3)}{\pi^2}
  \right).
\]

The product \(F_m(x)\) need not be used only at the zero \(x=1\).  Away from integral values
of \(x\), it gives convergent product values without removing any vanishing factor.  For
example, the first odd member gives
\[
  F_1(x) = \prod_{n=1}^{\infty} e^{x^2/n} \left(1-\frac{x^2}{n^2}\right)^n .
\]
At the half-value \(x=1/2\), this yields
\[
\prod_{n=1}^{\infty}
e^{1/(4n)}
\left(1-\frac{1}{4n^2}\right)^n
=
2^{1/12}A^{-3}
\exp\left(\frac{\gamma+2}{4}\right),
\]
where \(A\) is the Glaisher--Kinkelin constant.  This suggests that rational evaluations of
the same \(x\)-dependent Wallis products form a separate direction, distinct from the
products \(P_m\) considered here.

The two-factor extension also gives useful specialisations.  The simplest new case is the
pure hyperbolic companion, obtained by taking \(A=0\) and \(B=1\).  For \(m=0\),
\[
  P_0(0,1) = \prod_{n=2}^{\infty} \left(1+\frac1{n^2}\right) = \frac{\sinh\pi}{2\pi}.
\]
Thus the symmetric product with \(A=B=1\) gives
\[
  P_0(1,1) = \prod_{n=2}^{\infty} \left(1-\frac1{n^4}\right) = \frac{\sinh\pi}{4\pi}.
\]

For the first non-trivial even case, put
\[
  C_3 = \zeta(3) - \sum_{s=0}^{2} \frac{(2\pi)^s}{s!} \operatorname{Li}_{3-s}(e^{-2\pi}).
\]
Then the pure hyperbolic companion is
\[
P_2(0,1)
=
\prod_{n=2}^{\infty}
e^{-1}
\left(1+\frac1{n^2}\right)^{n^2}
=
\exp\left(
\frac32
-
\log2
-
\frac{\pi}{3}
-
\frac{C_3}{2\pi^2}
\right).
\]
Combining this with the ordinary \(P_2\) product gives the symmetric two-factor case
\[
\begin{aligned}
  P_2(1,1)
  &=
  \prod_{n=2}^{\infty}
  \left(1-\frac1{n^4}\right)^{n^2}
  \\
  &=
  \frac{\pi}{e^{3/2}}
  \exp\left(
  \frac32
  -
  \log2
  -
  \frac{\pi}{3}
  -
  \frac{C_3}{2\pi^2}
  \right)
  \\
  &=
  \frac{\pi}{2}
  \exp\left(
  -
  \frac{\pi}{3}
  -
  \frac{C_3}{2\pi^2}
  \right).
\end{aligned}
\]
The quotient case \(A=1\), \(B=-1\) gives
\[
\begin{aligned}
  P_2(1,-1)
  &=
  \prod_{n=2}^{\infty}
  e^2
  \left(
  \frac{1-\frac1{n^2}}{1+\frac1{n^2}}
  \right)^{n^2}
  \\
  &=
  \pi
  \exp\left(
  -
  3
  +
  \log2
  +
  \frac{\pi}{3}
  +
  \frac{C_3}{2\pi^2}
  \right).
\end{aligned}
\]

For the next even case, put
\[
  C_5 = \zeta(5) - \sum_{s=0}^{4} \frac{(2\pi)^s}{s!} \operatorname{Li}_{5-s}(e^{-2\pi}).
\]
Then
\[
P_4(0,1)
=
\prod_{n=2}^{\infty}
\exp\left(-n^2+\frac12\right)
\left(1+\frac1{n^2}\right)^{n^4}
=
\exp\left(
\frac14
-
\log2
+
\frac{\pi}{5}
+
\frac{3C_5}{2\pi^4}
\right).
\]
The symmetric two-factor case is
\[
\begin{aligned}
  P_4(1,1)
  &=
  \prod_{n=2}^{\infty}
  e
  \left(1-\frac1{n^4}\right)^{n^4}
  \\
  &=
  \pi
  \exp\left(
  -\frac32
  -
  \log2
  +
  \frac{3\zeta(3)}{\pi^2}
  +
  \frac{\pi}{5}
  +
  \frac{3C_5}{2\pi^4}
  \right).
\end{aligned}
\]


\begin{thebibliography}{99}

  \bibitem{wallis}
J. Wallis, \emph{Arithmetica Infinitorum}, Oxford, 1656.

  \bibitem{holcombe}
S. R. Holcombe, ``A product representation of pi,'' \emph{American Mathematical Monthly}
\textbf{120}, no. 8 (2013), 705. arXiv:1204.2451v3 [math.NT].

  \bibitem{allouche}
J.-P. Allouche, ``Hölder and Kurokawa meet Borwein--Dykshoorn and Adamchik,'' \emph{Journal
of the Ramanujan Mathematical Society} \textbf{38}, no. 3 (2023), 265--273.
arXiv:2205.09492v1 [math.NT].

  \bibitem{caihukim}
Q. Cai, S. Hu and M.-S. Kim, ``Euler's transformation, zeta functions and generalizations of
Wallis' formula,'' arXiv:2201.09674v4 [math.NT], 2022.

  \bibitem{farrell}
J. W. E. Farrell, ``Generalising the Wallis Product,'' arXiv:1906.00122 [math.GM], 2019.

  \bibitem{kurokawawakayama1}
N. Kurokawa and M. Wakayama, ``Duplication formulas in triple trigonometry,''
\emph{Proceedings of the Japan Academy, Series A, Mathematical Sciences} \textbf{79}, no. 8
(2003), 123--127.

  \bibitem{kurokawawakayama2}
N. Kurokawa and M. Wakayama, ``Extremal values of double and triple trigonometric
functions,'' \emph{Kyushu Journal of Mathematics} \textbf{58}, no. 1 (2004), 141--166.

  \bibitem{whittakerwatson}
E. T. Whittaker and G. Watson, \emph{A Course of Modern Analysis}, 4th ed., Cambridge
University Press, Cambridge, 1927.

\end{thebibliography}
\end{document}